\documentclass[11pt]{article}
\usepackage{mathrsfs}
\usepackage{epic,eepic,epsf,epsfig}
\usepackage{amsfonts,srcltx,mathrsfs}
\usepackage{amsmath}
\usepackage{amssymb}
\usepackage{amsbsy}
\usepackage{float}
\usepackage{graphicx}
\usepackage{amsfonts}
\usepackage{color}

\usepackage{url}
\newtheorem{lem}{Lemma}[section]
\newtheorem{theorem}[lem]{Theorem}

\newtheorem{case}{Case}
\newtheorem{prop}[lem]{Proposition}

\def\a{\alpha} \def\b{\beta}   
\def\r{\rho} \def\s{\sigma}

 \def\lg{\langle} \def\rg{\rangle}
\def\nd{\mathrel{\bigm|\kern-.7em/}}

\def\f{\noindent}

\def\Aut{\hbox{\rm Aut}}

\def\Aut{\hbox{\rm Aut}}

\newcommand{\qed}{\mbox{\raisebox{0.7ex}{\fbox{}}} \vspace{4truemm}}

\begin{document}

\title{On regular polytopes of order $2^n$}

\author{ \\ Dong-Dong Hou$^{a}$, Yan-Quan Feng*$^{a}$, Dimitri Leemans$^{b}$\\
$^{a}${\small Department of Mathematics, Beijing Jiaotong University, Beijing,
100044, P.R. China}\\
$^{b}${\small {D\'epartement de Math\'ematique, Universit\'e Libre de Bruxelles, 1050 Bruxelles Belgium}}}

\date{}
\maketitle

\footnotetext{*Corresponding author. \\ \mbox{}\hskip 0.8cm E-mails: yqfeng$@$bjtu.edu.cn, holderhandsome$@$bjtu.edu.cn,dleemans$@$ulb.ac.be
}
\begin{abstract}
For each $d\geq 3$, $n \geq 10$, and $k_1, k_2, \ldots, k_{d-1}\geq 2$ with $k_1+k_2+\ldots+k_{d-1}\leq n-1$, we construct a regular $d$-polytope whose automorphism group is of order $2^n$ and whose Schl\"afli type is $\{2^{k_1},2^{k_2}, \ldots, 2^{k_{d-1}}\}$.

\bigskip
\f {\bf Keywords:} Regular polytope, $2$-group, automorphism group, string C-group.
\end{abstract}

\section{Introduction}

In the 1960s, Branko Gr\"unbaum suggested to the geometric community to study generalizations of the concept of regular polytopes that he called polystromata. His work greatly influenced Ludwig Danzer and Egon Schulte who developed, along those lines, the theory of what are now called abstract regular polytopes. The comprehensive book written by Peter McMullen and Egon Schulte~\cite{ARP} is nowadays seen as the reference on the subject.

Abstract polytopes are a generalization of the classical notion of convex geometric polytopes to more general structures. The highly symmetric examples are the most studied. They
include not only classical regular polytopes such as the Platonic solids, but also non-degenerate
regular maps on surfaces. Another famous example is the 11-cell that Gr\"unbaum discovered in 1977~\cite{grunbaum} by gluing together eleven hemi-icosahedra in such a way that the ``geometry around each vertex" would look like a hemi-dodecahedron. 

An abstract regular polytope is a partially-ordered set endowed with a rank function, satisfying certain conditions that arise naturally from a geometric setting. 
There are numerous references on regular polytopes in the literature. The order of the automorphism group of a regular polytope is also called the {\em order} of the regular polytope. The atlas~\cite{atles1} contains information about all regular polytopes with order at most 2000. Up to now, most constructions of regular polytopes were obtained either by computer searches or by using almost simple groups as their automorphism groups. We refer to \cite{psl3q,psl4q,CFLM2017,soclepsl2q,fl,flm1,sympolcorr,flm2,flm,suzuki,psl2q,pgl2q} for examples.
Most of the theoretical results mentioned above were inspired by experimental data collected over the years in atlases of polytopes such as~\cite{atles1}. Among them, the authors would like to single out reference~\cite{Lalg} which was accepted for publication by Gr\"unbaum, who was a strong supporter of that approach of collecting experimental data, analyzing them and stating conjectures that could then be proved by developing new mathematical tools.

There are  just a few theoretical constructions of regular polytopes for solvable groups (see \cite{sc1024,HFL,Loyola}). In~\cite{Problem} Schulte and Weiss proposed the following problem: for a positive integer $n$, characterize regular polytopes of orders $2^n$.

Let $d\geq 3$, $n \geq 10$, and $k_1, k_2, \ldots, k_{d-1}\geq 2$. Let $\mathcal{P}$ is a regular $d$-polytope of order $2^n$ and type $\{2^{k_1},2^{k_2}, \ldots, 2^{k_{d-1}}\}$. We know, from Conder~\cite{SmallestPolytopes} (see Proposition~\ref{leastorder} below), that $k_1+k_2+\ldots+k_{d-1}\leq n-1$. Cunningham and Pellicer~\cite{GD2016} classified regular $3$-polytopes with $k_1+k_2=n-1$, and the authors constructed in~\cite{HFL} a regular $3$-polytope for each $k_1,k_2, n$ with $k_1+k_2\leq n-1$.
In this paper, for each $d\geq 4$, we construct a regular $d$-polytope for each $k_1, k_2, \ldots, k_{d-1},n$ with $k_1+k_2+\ldots+k_{d-1}\leq n-1$. 
Our main theorem can be stated as follows.

\begin{theorem}\label{existmaintheorem}
For any integers $d,n, k_1, k_2, \ldots, k_{d-1}$ such that $d\geq 3$, $n \geq 10$, $k_1, k_2, \ldots,$
$ k_{d-1}\geq 2$ and $k_1+k_2+\ldots +k_{d-1}\leq n-1$,
there exists a string C-group
$(G, \{\r_0, \r_1, \ldots, \r_{d-1}\})$ of order $2^n$ and type $\{2^{k_1}, 2^{k_2}, \ldots, 2^{k_{d-1}}\}$.
\end{theorem}

The paper is organized as follows. In Section~\ref{backgroud}, we give the necessary definitions to understand this paper and we recall some results that we use in Section~\ref{Main Results} to prove Theorem~\ref{existmaintheorem}.

\section{Background definitions and preliminaries}\label{backgroud}

Abstract regular polytopes and string C-groups are the same mathematical objects. The link between these objects may be found for instance in~\cite[Section 2E]{ARP}.
We take here the viewpoint of string C-groups because it is the easiest and the most efficient one to define abstract regular polytopes.
Let $G$ be a group and let $S=\{\rho_0,\ldots,\rho_{d-1}\}$ be a generating set of involutions of $G$.
For $I\subseteq\{0,\ldots,d-1\}$, let $G_I$ denote the group generated by $\{\rho_i:i\in I\}$.
Suppose that
\begin{itemize}
\item[*]  for any $i,j\in \{0, \ldots, d-1\}$ with $|i-j|>1$, $\rho_i$ and $\rho_j$ commute (the \emph{string
property});
\item[*] for any $I,J\subseteq\{0,\ldots,d-1\}$,
$G_I\cap G_J=G_{I\cap J}\ \  (\mbox{the \emph{intersection property}})$.
\end{itemize}
Then the pair $(G,S)$ is called a {\em string C-group of rank} $d$ and the {\em order} of $(G,S)$ is simply the order of $G$. The {\em type} of $(G,S)$ is the ordered set $\{k_1,\ldots, k_{d-1}\}$, where $k_i$ is the order of $\r_{i-1}\r_i$. It is natural to assume that each $k_i$ is at least 3 for otherwise the generated group is a direct product of two smaller groups and the corresponding string C-group is called {\em degenerate}. By the intersection property, $S$ is a minimal generating set of $G$.

If $(G,S)$  only satisfies the string property, it is called a {\em string group generated by involutions}, or an \emph {sggi} for short. The following proposition is called the {\em quotient criterion} for a string C-group.

\begin{prop}{\rm \cite[Theorem 2E17]{ARP}}\label{stringC}
Let $\Gamma=\lg \r_0, \r_1, \ldots, \r_{d-1}\rg$ be an sggi, and let $(\Delta,$ $\{\s_0, \s_1,\ldots, \s_{d-1}\})$ be a string C-group.
If the mapping $\r_j \mapsto \s_j$ for $j=0, \ldots, d-1$ induces a homomorphism $\pi : \Gamma \rightarrow \Delta$, which is one-to-one on the subgroup
$\Gamma_{d-1}=\lg \r_0, \r_1, \ldots, \r_{d-2}\rg$ or on $\Gamma_{0}=\lg \r_1, \r_2, \ldots, \r_{d-1}\rg$, then $(\Gamma, \{\r_0, \r_1, \ldots, \r_{d-1}\})$
is also a string C-group.
\end{prop}

Let $(G,S)$ be a string C-group and let $\mathcal{P}$ be its corresponding regular polytope. Then the rank, the order, and the type of $(G,S)$ mean the {\em rank}, the {\em order}, and the {\em (Schl\"afli) type} of $\mathcal{P}$, respectively. A regular polytope $\mathcal{P}$ is called a regular {\em $d$-polytope}, if  $\mathcal{P}$ has rank $d$. Conder~\cite{SmallestPolytopes} obtained a lower bound for the order of a regular polytope.

\begin{prop}{\rm \cite[Theorem 3.2]{SmallestPolytopes}}\label{lowerbound}\label{leastorder}
If $\mathcal{P}$ is a regular $d$-polytope of type $\{k_1, k_2, \ldots, k_{d-1}\}$, then $\mathcal{P}$ has order at least $2k_1k_2 \ldots k_{d-1}$.
\end{prop}

If the lower bound in Proposition~\ref{lowerbound} is attained, $\mathcal{P}$ is called {\em tight}. A {\em string Coxeter group $[k_1, k_2, \ldots, k_{d-1}]$} is defined as the following group:
\begin{eqnarray*}
  & \lg \r_0, \r_1, \ldots, \r_{d-1} \ |\ & \r_i^2=1 \ {\rm for} \ 0 \leq i \leq d-1, (\r_i\r_{i+1})^{k_{i+1}}=1 \ {\rm for} \ 0 \leq i \leq d-2,  \\
  &  & (\r_i\r_j)^2=1 \ {\rm for} \ 0 \leq i < j-1  < d-1\rg.
\end{eqnarray*}

\begin{prop}{\rm \cite[Theorem 5.3]{SmallestPolytopes}}\label{tight}
For every sequence $(k_1, k_2, \ldots, k_{d-1})$ of $d-1$ even integers greater than $2$, there exists a tight regular $d$-polytope $\mathcal{P}$ of order $2k_1k_2 \ldots k_{d-1}$ and type $\{k_1, k_2, \ldots, k_{d-1}\}$. In particular, one can take the string Coxeter group $[k_1, k_2, \ldots, k_{d-1}]$
with standard generators $\r_0, \r_1, \ldots, \r_{d-1}$, and let $\Aut(\mathcal{P})$ be the quotient obtained by adding all relations of the form $[\r_{i}, (\r_{i+1}\r_{i+2})^2]=1$
for $0 \leq i \leq d-3$ and $[(\r_i\r_{i+1})^2, \r_{i+2}]$ for $0 \leq i \leq d-3$.
\end{prop}

The following proposition gives some constructions for string C-groups of order $2^n$.

\begin{prop}{\rm \cite[Theorem 1.2]{HFL}}\label{$2^s,2^t$}
Let $n \geq 10$, $s, t \geq 2$ and $n-s-t \geq 1$. Set $R= \{\r_0^2, \r_1^2, \r_2^2, (\r_0\r_1)^{2^{s}}, (\r_1\r_2)^{2^{t}}, (\r_0\r_2)^2$, $[(\r_0\r_1)^4, \r_2], [\r_0,(\r_1\r_2)^4]\}$ and define
$$H=\left\{
\begin{array}{ll}
\lg \r_0, \r_1, \r_2 \ |\ R, [(\r_0\r_1)^2, \r_2]^{2^{\frac{n-s-t-1}{2}}}\rg, & n-s-t\mbox{ odd }\\
\lg \r_0, \r_1, \r_2 \ |\ R, [(\r_0\r_1)^2, (\r_1\r_2)^2]^{2^{\frac{n-s-t-2}{2}}}\rg, & n-s-t\mbox{ even. }
\end{array}
\right.$$
Then $(H,\{\r_0,\r_1,\r_2\})$ is a string C-group of order $2^n$ and type $\{ 2^s,2^t\}$.
\end{prop}

Let $G$ be a group. For $a, b\in G$, we use $[a, b]$ as an abbreviation for the
{\em commutator} $a^{-1}b^{-1}ab$ of $a$ and $b$. The following  proposition is a basic property of commutators and its proof is straightforward.

\begin{prop}\label{commutator}
Let $G$ be a group. Then, for any $a, b, c \in G$, $[ab, c]=[a, c]^b[b, c]$ and $[a, bc]=[a, c][a, b]^c$.
\end{prop}

Finally, we will also use the following result in the proof of our theorem.

\begin{prop}\label{claim}
Let $H=\lg a,b,c\rg$ be a group such that $a^2=b^2=c^2=(ac)^2=[(ab)^2, c]=1$. Then $[a, (bc)^2]=[a, (bc)^4]=[(ab)^4, c]=1$, $\lg (ab)^2 \rg \unlhd H$ and $\lg (bc)^2 \rg\unlhd H$.
\end{prop}
{\bf Proof:} Since $(ac)^2=1$, we have $[a, c]=1$. By Proposition~\ref{commutator}, $1=[(ab)^2, c]=[a, c]^{bab}[bab, c]=[bab, c]=[a, bcb]^b=[a, bcbcc]^b=[a, c]^b[a, (bc)^2]^{cb}=[a, (bc)^2]^{cb}$, that is, $[a, (bc)^2]=1$. It follows that $[(ab)^4, c]=1$ and $[a, (bc)^4]=1$.
Note that $(ab)^a=ba=(ab)^{-1}$ and $(ab)^b=ba=(ab)^{-1}$. Since $H=\lg a,b,c\rg$, we have $\lg (ab)^2 \rg \unlhd H$, and similarly, $\lg (bc)^2 \rg \unlhd H$, as required.
\hfill\qed

\section{Proof of Theorem~\ref{existmaintheorem}}\label{Main Results}

The proof of this theorem is constructive.
Let $d,n, k_1, k_2, \ldots, k_{d-1}$ be integers such that $d\geq 3$, $n \geq 10$, and $k_1, k_2, \ldots, k_{d-1}\geq 2$ with $k_1+k_2+\ldots +k_{d-1}\leq n-1$. Define

$$G=\left\{
\begin{array}{ll}
\lg \r_0, \r_1, \ldots, \r_{d-1} \ |\ R_1, R_2, R_3, [(\r_{d-3}\r_{d-2})^2, \r_{d-1}]^{2^{\frac{l-1}{2}}}\rg, & \mbox{\ \ \  $l$ odd }\\
\lg \r_0, \r_1, \ldots, \r_{d-1} \ |\ R_1, R_2, R_3, [(\r_{d-3}\r_{d-2})^2, (\r_{d-2}\r_{d-1})^2]^{2^{\frac{l-2}{2}}}\rg, & \mbox{\ \ \ $l$ even,}
\end{array}
\right.$$

\f where
\begin{eqnarray}
    R_1=\{\r_i^2, (\r_j\r_k)^2, (\r_\ell\r_{\ell+1})^{2^{k_{\ell+1}}}  \ |\ 0 \leq i \leq  d-1,  0 \leq j < k-1  \leq d-2, 0\leq \ell\leq d-2\}, && \label{eq1} \\
  R_2=\{[(\r_i\r_{i+1})^2, \r_{i+2}] \ |\ 0 \leq i \leq  d-4\}\}, \label{eq2} \hskip 8.1cm &&\\  R_3=\{[\r_{d-3},
   (\r_{d-2}\r_{d-1})^4], [(\r_{d-3}\r_{d-2})^4, \r_{d-1}]\},\hskip 7.2cm && \label{eq3}
\end{eqnarray}
and $l=n-(k_1+k_2+\ldots +k_{d-1})$. Note that $l\geq 1$.

The case where $d=3$ has been dealt with in Proposition~\ref{$2^s,2^t$}. Hence we may assume from now on that $d\geq 4$. 
For convenience, write $o(g)$ for the order of $g$ in $G$. 

Let $$K=\lg \r_0, \r_1, \ldots, \r_{d-1} \ |\ R_1, R_2, [(\r_{d-3}\r_{d-2})^2, \r_{d-1}]\rg,$$
where $R_1$ and $R_2$ are given by equations (\ref{eq1}) and (\ref{eq2}). Then $[(\r_{i}\r_{i+1})^2, \r_{i+2}]=1$ for
$0 \leq i \leq d-3$ in $K$, and by Proposition~\ref{claim}, $[\r_i, (\r_{i+1}\r_{i+2})^2]=1$. Then $(K, \{\r_0, \r_1, \ldots, \r_{d-1}\})$
is a string C-group of order $2^{1+k_1+k_2+\ldots+k_{d-1}}$ by Proposition~\ref{tight}, with type $\{2^{k_1}, 2^{k_2}, \ldots, 2^{k_{d-1}}\}$. In particular, the listed exponents are the true orders of the corresponding elements in $K$, that is, $o(\r_i)=2$ for $0 \leq i \leq d-1$, $o(\r_i\r_{i+1})=2^{k_{i+1}}$ for $0 \leq i \leq d-2$ and $o(\r_i\r_j)=2$ for $0 \leq i < j-1  < d-1$.

Note that $[(\r_{d-3}\r_{d-2})^2, \r_{d-1}]=1$ in $K$. Proposition~\ref{claim} implies that $[(\r_{d-3}\r_{d-2})^4, \r_{d-1}]=[\r_{d-3}, (\r_{d-2}\r_{d-1})^4]=1$ in $K$. Furthermore, $[(\r_{d-3}\r_{d-2})^2, (\r_{d-2}\r_{d-1})^2]=[(\r_{d-3}\r_{d-2})^2, \r_{d-1}]\newline [(\r_{d-3}\r_{d-2})^2, \r_{d-2}\r_{d-1}\r_{d-2}]^{\r_{d-1}}=[(\r_{d-3}\r_{d-2})^{-2}, \r_{d-1}]^{\r_{d-2}\r_{d-1}}=1$. In particular, we have, in $K$, that $[(\r_{d-3}\r_{d-2})^2, \r_{d-1}]^{2^{\frac{l-1}{2}}}=1$ when $l$ is odd, and $[(\r_{d-3}\r_{d-2})^2, (\r_{d-2}\r_{d-1})^2]^{2^{\frac{l-2}{2}}}=1$ when $l$ is even. Thus $\r_{0}, \r_{1}, \ldots, \r_{d-1}$ in $K$ satisfy the same relations as do $\r_{0}, \r_{1}, \ldots, \r_{d-1}$ in $G$, and hence the map: $\r_i \mapsto \r_i$, $0\leq i \leq d-1$, induces an epimorphism $\alpha$ from $G$ to $K$. This also implies, in $G$, that $o(\r_i)=2$ for $0 \leq i \leq d-1$, $o(\r_i\r_{i+1})=2^{k_{i+1}}$ for $0 \leq i \leq d-2$, and $o(\r_i\r_j)=2$ for $0 \leq i < j-1  < d-1$, because these are true in $K$.

Let $K_1=\lg \r_0, \r_1, \ldots, \r_{d-2}\rg\leq K$ and $G_1=\lg \r_0, \r_1, \ldots, \r_{d-2}\rg\leq G$. Then the restriction $\a |_{G_1}$ of $\a$ on $G_1$ is an epimorphism from $G_1$ to $K_1$.
Now we prove that $\a |_{G_1}$ is actually an isomorphism from $G_1$ to $K_1$. To do that, let
$$L=\lg \r_0, \r_1, \ldots, \r_{d-2}\ |\ \overline{R}_1,R_2 \rg, $$
where $\overline{R}_1=\{\r_i^2, (\r_j\r_k)^2, (\r_\ell\r_{\ell+1})^{2^{k_{\ell+1}}}  \ |\ 0 \leq i \leq  d-2,  0 \leq j < k-1  \leq d-3, 0\leq \ell\leq d-3\}$ and $R_2$ is given by equation (\ref{eq2}).
Now, a easy argument similar to $K$ shows that $|L|=2^{1+k_1+k_2+\ldots+k_{d-2}}$.

Clearly, $\r_0, \r_1, \ldots, \r_{d-2},1$ in $L$ satisfy the same relations as do $\r_0, \r_1, \ldots, \r_{d-2},\r_{d-1}$ in $K$, and therefore, the map $\r_i\mapsto \r_i$ for $0\leq i\leq d-2$, with $\r_{d-1}\mapsto 1$, induces an epimorphism from $K$ to $L$, whose restriction on $K_1$ is an epimorphism from $K_1$ to $L$. On the other hand, $\r_0, \r_1, \ldots, \r_{d-2}$ in $K_1\leq K$ satisfy the same relations as do $\r_0, \r_1, \ldots, \r_{d-2}$ in $L$, and hence there is an epimorphism from $L$ to $K_1$.  These facts yield that the map $\r_i\mapsto \r_i$ for $0\leq i\leq d-2$ induces an isomorphism from $K_1$ to $L$, because $|L|$ is finite.

Similarly, $\r_0, \r_1, \ldots, \r_{d-2},1$ in $L$ satisfy the same relations as do $\r_0, \r_1, \ldots, \r_{d-2},\r_{d-1}$ in $G$, and $\r_0, \r_1, \ldots, \r_{d-2}$ in $G_1 \leq G$ satisfy the same relations as do $\r_0, \r_1, \ldots, \r_{d-2}$ in $L$. A similar argument to the one of the above paragraph permits to conclude that the map $\r_i\mapsto \r_i$ for $0\leq i\leq d-2$ induces an isomorphism from $G_1$ to $L$. This, together with the isomorphism from $K_1$ to $L$ in the above paragraph, implies that $\a |_{G_1}$ is an isomorphism from $G_1$ to $K_1$.

Since $(K, \{\r_0, \r_1, \ldots, \r_{d-1}\})$ is a string C-group, $(G, \{\r_0, \r_1, \ldots, \r_{d-1}\})$ is a string C-group by Proposition~\ref{stringC}, which has type $\{2^{k_1}, 2^{k_2}, \ldots, 2^{k_{d-1}}\}$. To finish the proof, we are only left to show that $|G|=2^n$. We prove this by induction on $d$. It is true for $d=3$ by Proposition~\ref{$2^s,2^t$}, and we may let $d\geq 4$.

Let $N=\lg (\r_0\r_1)^2\rg \leq G$. Since $o(\r_0\r_1)=2^{k_1}$ in $G$, we have $|N|=2^{k_1-1}$. By Proposition~\ref{claim}, $N\unlhd G$, because $[\r_0,\r_j]=[\r_1,\r_j]=1$ for any $j\geq 3$.
Clearly, $G/N\cong M$ with
$$M=\left\{
\begin{array}{ll}
\lg \r_0, \r_1, \ldots, \r_{d-1} \ |\ R_1, R_2, R_3, [(\r_{d-3}\r_{d-2})^2, \r_{d-1}]^{2^{\frac{l-1}{2}}},(\r_0\r_1)^2\rg, & \mbox{\ \ \  $l$ odd }\\
\lg \r_0, \r_1, \ldots, \r_{d-1} \ |\ R_1, R_2, R_3, [(\r_{d-3}\r_{d-2})^2, (\r_{d-2}\r_{d-1})^2]^{2^{\frac{l-2}{2}}},(\r_0\r_1)^2\rg, & \mbox{\ \ \ $l$ even,}
\end{array}
\right.$$
where $R_1$, $R_2$, $R_3$ are given by equations~(\ref{eq1}),~(\ref{eq2}) and (\ref{eq3}). Write $M_1=\lg \r_1,\r_2,\ldots,\r_{d-1}\rg\leq M$.

Since  $o(\r_0\r_1)=2^{k_1}\geq 4$ in $G$, $\lg \r_0,\r_1\rg$ is a dihedral group of order $2^{k_1+1}\geq 8$, implying that $\r_0\not\in N$ and $\r_1\not\in N$. In particular, $o(\r_0N)=2$ in $G/N$, and therefore,  $o(\r_0)=2$ in $M$. Since $(G, \{\r_0, \r_1, \ldots, \r_{d-1}\})$ is a string C-group, $\lg \r_1,\r_2,\ldots,\r_{d-1}\rg\cap \lg \r_0\rg N\leq \lg \r_1,\r_2,\ldots,\r_{d-1}\rg\cap \lg \r_0,\r_1\rg=\lg \r_1\rg$. If $\r_1\in  \lg \r_0\rg N$ then either $\r_1\in N$ or $\r_0\r_1\in N$, both of which are impossible. It follows that $\lg \r_1,\r_2,\ldots,\r_{d-1}\rg\cap \lg \r_0\rg N=1$ in $G$, and hence $\lg \r_1,\r_2,\ldots,\r_{d-1}\rg N\cap \lg \r_0\rg N=N(\lg \r_1,\r_2,\ldots,\r_{d-1}\rg\cap \lg \r_0\rg N)=N$ in $G/N$. This implies that $\lg \r_0\rg\cap M_1=1$ in $M$. Since $(\r_0\r_j)^2=1$ in $M$ for any $j\geq 2$, we have $M=\lg \r_0\rg\times M_1$ and hence $|M|=2|M_1|$.

Let
$$A=\left\{
\begin{array}{ll}
\lg \r_1, \ldots, \r_{d-1} \ |\ R_1^-, R_2^-, R_3, [(\r_{d-3}\r_{d-2})^2, \r_{d-1}]^{2^{\frac{l-1}{2}}}\rg, & \mbox{\ \ \  $l$ odd }\\
\lg \r_1, \ldots, \r_{d-1} \ |\ R_1^-, R_2^-, R_3, [(\r_{d-3}\r_{d-2})^2, (\r_{d-2}\r_{d-1})^2]^{2^{\frac{l-2}{2}}}\rg, & \mbox{\ \ \ $l$ even,}
\end{array}
\right.$$
\f where $R_1^-=\{\r_i^2, (\r_j\r_k)^2, (\r_\ell\r_{\ell+1})^{2^{k_{\ell+1}}}  \ |\ 1 \leq i \leq  d-1,  1 \leq j < k-1  \leq d-2, 1\leq \ell\leq d-2\}$,
$R_2^-=\{[(\r_i\r_{i+1})^2, \r_{i+2}] \ |\ 1 \leq i \leq  d-4\}\}$ and $R_3$ is given by equation~(\ref{eq3}). Note that $l=(n-k_1)-k_2-\ldots-k_{d-1}$.

Let $n-k_1\geq 10$. Since $(G, \{\r_0, \r_1, \ldots, \r_{d-1}\})$ is a string C-group, by taking $n-k_1$ in $A$ as $n$ in $G$, $(A, \{\r_1, \r_2, \ldots, \r_{d-1}\})$ is a string C-group of rand $d-1$. Then the inductive hypothesis implies that $|A|=2^{n-k_1}$. Now we claim that this is also true for $n-k_1\leq 9$.

Note that $d-1\geq 3$, $l \geq 1$ and $k_1, k_2, \ldots, k_{d-1} \geq 2$. Since $1+2(d-2)\leq l+k_2+k_3+\ldots +k_{d-1}=n-k_1 \leq 9$, we have $3\leq d-1 \leq 5$ and $5\leq n-k_1\leq 9$.

First assume $d-1=3$. Then $A=\lg \r_1, \r_2, \r_3\rg$ and $l=(n-k_1)-k_2-k_3$. If $n-k_1=5$, then $(l, k_2, k_3)=(1,2,2)$, and using {\sc Magma}~\cite{BCP97} we easily check that $|A|=2^5=2^{n-k_1}$.
If $n-k_1=6$, then $(l, k_2, k_3)=(1,2,3)$, $(1,3,2)$, or $(2,2,2)$; if $n-k_1=7$, then $(l, k_2, k_3)=(1,2,4)$, $(1,4,2)$, $(1,3,3)$, $(2,2,3)$, $(2,3,2)$, or $(3,2,2)$; if $n-k_1=8$, then $(l, k_2, k_3)=(1,2,5)$, $(1,5,2)$, $(1,3,4)$, $(1,4,3)$, $(2,2,4)$, $(2,4,2)$, $(2,3,3)$, $(3,2,3)$, $(3,3,2)$ or $(4,2,2)$; if $n-k_1=9$, then $(l, k_2, k_3)=(1,2,6)$, $(1,6,2)$, $(1,3,5)$, $(1,5,3)$, $(1,4,4)$, $(2,2,5)$, $(2,5,2)$,
$(2,3,4)$, $(2,4,3)$, $(3,2,4)$, $(3,4,2)$, $(3,3,3)$, $(4,2,3)$, $(4,3,2)$, or $(5,2,2)$. For each $(l, k_2, k_3)$, {\sc Magma} computations show that  $|A|=2^{n-k_1}$.

Assume $d-1=4$. Then $A=\lg \r_1, \r_2, \r_3, \r_4\rg$ and $7\leq (n-k_1)=l+k_2+k_3+k_4\leq 9$. If $n-k_1=7$, then  $(l, k_2, k_3, k_4)=(1,2,2,2)$; if $n-k_1=8$, then $(l, k_2, k_3, k_4)=(1,3,2,2)$, $(1,2,3,2)$, $(1,2,2,3)$ or $(2,2,2,2)$; if $n-k_1=9$, then $(l, k_2, k_3, k_4)=(1,4,2,2)$, $(1,2,4,2)$, $(1,2,2,4)$, $(1,3,3,2)$, $(1,3,2,3)$, $(1,2,3,3)$, $(2,3,2,2)$, $(2,2,3,2)$, $(2,2,2,3)$ or $(3,2,2,2)$. Assume $d-1=5$. Then $A=\lg \r_1, \r_2, \r_3, \r_4, \r_5\rg$ and $n-k_1=9$; furthermore, $(l, k_2, k_3, k_4, k_5)=(1,2,2,2,2)$. Again using {\sc Magma}~\cite{BCP97}, for each case we have $|A|=2^{n-k_1}$, as claimed.

Clearly, $1,\r_1,\r_2,\ldots,\r_{d-1}$ in $A$ satisfy the same relations as $\r_0,\r_1,\r_2,\ldots,\r_{d-1}$ in $M$. Thus the map $\r_0\mapsto 1$, $\r_i\mapsto \r_i$ for $1\leq i\leq d-1$, induces an epimorphism $\beta$ from $M$ to $A$ and hence the restriction $\b|_{M_1}$ is an epimorphism from $M_1$ to $A$. On the other hand, $\r_1,\r_2,\ldots,\r_{d-1}$ in $M_1 \leq M$ satisfy the same relations as $\r_1,\r_2,\ldots,\r_{d-1}$ in $A$, and therefore, there is an epimorphism from $A$ to $M_1$. Thus, $\b|_{M_1}$ is an isomorphism from $M_1$ to $A$ and in particular, $|M_1|=|A|=2^{n-k_1}$. Now, we have $|G|=|G/N||N|=|M||N|=2|M_1||N|=2\cdot 2^{n-k_1}\cdot 2^{k_1-1}=2^n$. This completes the proof.\hfill\qed

\f {\bf Acknowledgements:} This work was supported by the National Natural Science Foundation of China (11571035, 11731002) and the 111 Project of China (B16002).

\end{document}